\documentclass[reprint,aps,pre,floatfix]{revtex4-1}

\usepackage{xcolor}
\usepackage[normalem]{ulem}
\usepackage{amsmath,amsfonts,amscd,amssymb,graphicx}
\usepackage{lineno}
\usepackage{dcolumn}% Align table columns on decimal point
\usepackage{overpic}
\usepackage{color}
\usepackage{graphicx}
\usepackage{mathrsfs}
\usepackage{subcaption}
\usepackage{adjustbox}
\usepackage{rotating}
\usepackage{float}
\makeatletter
\let\newfloat\newfloat@ltx
\makeatother
\usepackage{algorithm}
\usepackage{algpseudocode}
\algnewcommand\algorithmicinput{\textbf{Input:}}
\algnewcommand\Input{\item[\algorithmicinput]}

\begin{document}

\title{Unsupervised multi-scale diagnostics}

\author{Karl Lapo$^{1}$, Sara M. Ichinaga$^{2}$, and J. Nathan Kutz$^{3}$}
 \affiliation{$^{1}$Department of Atmospheric and Cryospheric Sciences, University of Innsbruck, Innsbruck, Austria, 6020} 
 \affiliation{$^{2}$ Department of Applied Mathematics, University of Washington, Seattle, WA 98195}  
 \affiliation{$^{3}$ Departments of Applied Mathematics and Electrical and Computer Engineering, University of Washington, Seattle, WA} 

\begin{abstract}
The unsupervised and principled diagnosis of multi-scale data is a fundamental obstacle in modern scientific problems from, for instance, weather and climate prediction, neurology, epidemiology, and turbulence. Multi-scale data is characterized by a combination of processes acting along multiple dimensions simultaneously, spatiotemporal scales across orders of magnitude, non-stationarity, and/or invariances such as translation and rotation. Existing methods are not well-suited to multi-scale data, usually requiring supervised strategies such as human intervention, extensive tuning, or selection of ideal time periods. We present the multi-resolution Coherent Spatio-Temporal Scale Separation (mrCOSTS), a hierarchical and automated algorithm for the diagnosis of coherent patterns or modes in multi-scale data. mrCOSTS is a variant of Dynamic Mode Decomposition which decomposes data into bands of spatial patterns with shared time dynamics, thereby providing a robust method for analyzing multi-scale data. It requires no training but instead takes advantage of the hierarchical nature of multi-scale systems. We demonstrate mrCOSTS using complex multi-scale data sets that are canonically difficult to analyze: 1) climate patterns of sea surface temperature, 2) electrophysiological observations of neural signals of the motor cortex, and 3) horizontal wind in the mountain boundary layer. With mrCOSTS, we trivially retrieve complex dynamics that were previously difficult to resolve while additionally extracting hitherto unknown patterns of activity embedded in the dynamics, allowing for advancing the understanding of these fields of study. This method is an important advancement for addressing the multi-scale data which characterize many of the grand challenges in science and engineering.
\end{abstract}

\maketitle

\section{Introduction}\label{Sec:Intro}

From electromagnetism and quantum mechanics to fluid dynamics, the observed patterns in nature have been critical for hypothesis formulation and testing \citep{Kutz_DMD-textoobk_2016}. Data has enabled the discovery and construction of the physical laws of nature, as well as uniscale descriptions at specific levels of temporal and spatial resolution (e.g. Schr\"odinger equation for quantum mechanics at angstrom scale, Maxwell equations at hundreds of nanometers for visible light, etc). The problem in modern science, especially in grand challenge science and engineering settings \citep{omennGrandChallengesGreat2006} are related to deducing patterns in the complexities of multi-scale data, where orders of magnitude of temporal and spatial scales are jointly measured and observed. Thus, there is no unifying uniscale model description, rather dynamics are driven by different processes at, for instance, micro-, meso-, and macro-scales, all of which are coupled together in ways that are poorly understood. Instead, there is a general scale-separation problem in which we cannot even properly identify the patterns of activity at each scale in order to generate and test hypotheses about the underlying dynamics and their couplings.

A multi-scale diagnostic is therefor critical for advancing and testing hypotheses about the underlying physical laws and their couplings in the grand challenges in science and engineering such as weather and climate modeling, neuroscience, biological systems, finance, etc. all of which inherently exhibit high-dimensional, multi-scale dynamics \citep{omennGrandChallengesGreat2006}. The interplay between deductive and inductive reasoning for scientific theories has always required clear observations of a process in order to advance a successful theory. It is simpler to explain patterns once it is understood what you are looking for (e.g., supervised learning), but when the existence of a pattern is undetected and unknown, the underlying dynamic process remains a latent contribution and true understanding is incomplete. This is why unsupervised algorithms for finding patterns are of such exceptional value: they fundamentally can help to accelerate the process of scientific understanding by the discovery of patterns hidden to observers. Multi-scale dynamics are perhaps the most difficult from which to extract patterns, as the significant interactions across temporal and spatial scales prevent a clear observation of the underlying processes. These complexities have likely contributed to a significant lack of algorithmic development and analysis tools for multi-scale phenomenon. Demonstrated here is an algorithm capable of extracting hitherto unobserved patterns of activity in exceptionally challenging multi-scale problems.

Partially as a response to the challenges of multi-scale dynamics, machine learning has emerged as a leading method for analysis. However, machine learning has clear drawbacks in terms of computational costs \citep{szeEfficientProcessingDeep2017} and interpretability \citep{carvalhoMachineLearningInterpretability2019, linardatosExplainableAIReview2021}. Thus, there is a clear need for unsupervised methods which automate the extraction of patterns and diagnoses of multi-scale data. Data-driven model discovery has emerged as powerful paradigm for handling data with complex dynamics by attempting to discover missing dynamics from data directly instead of from first principles \citep{Brunton2019}. These approaches are powerful even in fields for which the first principles are known, as behavior at a larger scale tends to emerge from complex interactions at smaller scales, such as is the case for fluids \cite{Taira2017, Taira2020}, turbulence \cite{Rowley2009_dmd}, weather \citep{shahVeryLargeScaleMotionsAtmospheric2014, finniganBoundaryLayerFlowComplex2020,  honnertAtmosphericBoundaryLayer2020}, neurology \cite{einevollModellingAnalysisLocal2013, herrerasLocalFieldPotentials2016a}, and climate \cite{schurerSeparatingForcedChaotic2013, frankcombeSeparatingInternalVariability2015, Deser2016a, zhangPacificDecadalOscillation2018}, among many others. However, real multi-scale data sets remain a challenge for data-driven model discovery.

It is often stated that data are high-dimensional or multi-scale, but these definitions are in general used loosely. Here, we explicitly define multi-scale dynamics as those characterized by a combination of at least two of the following properties: being multivariate (i.e. acting along multiple dimensions simultaneously), containing process scales across orders of magnitude, being non-stationary, or having invariances such as translation and rotation. Additionally, these data often contain noise or uncertainty, i.e. from instrument error. We call data characterized by multi-scale dynamics ``multi-scale data". Data with these properties are precisely those which are least amenable to analysis since no method can work with all of the aforementioned properties of multi-scale data. For instance, multi-resolution analysis (MRA) can provide multi-scale information on the dynamics along a single dimension (e.g., time or along a spatial dimension) but not multiple dimensions simultaneously \cite{Kutz_DMD-textoobk_2016}. Modal analysis can reveal dominant spatial patterns \cite{Taira2017, Taira2020} but these methods cannot robustly handle mulitvariate data \cite[e.g.,][]{dommengetCautionaryNoteInterpretation2002}, non-stationary \cite[e.g.][]{chenVariantsDynamicMode2012}, and/or translating processes \cite{Kutz2016_dmd}. As such, diagnosing multi-scale processes, especially in an unsupervised approach with minimal tuning, is a fundamental methodological need for a diverse set of disciplines. Consequently, many of the grand challenges in science and engineering are, in essence, challenges of multi-scale data \citep{omennGrandChallengesGreat2006}.

To this end, we present multi-resolution Coherent Spatiotemporal Scale-separation (mrCOSTS). It is specifically tailored for robust, interpretable diagnoses of multi-scale dynamics; a hierarchy of scales are robustly identified in an unsupervised fashion with minimal tuning. In short, mrCOSTS operates on a principle similar to the continuous wavelet transforms. A sliding window is applied to the data and in each window a Dynamic Mode Decomposition (DMD) model is fit (Part 1, Fig. \ref{fig:flow-chart}). DMD has the desirable property of describing coherent spatio-temporal modes. Each mode is interpreted as a coherent spatial mode governed by a single set of temporal dynamics and a small number of coherent modes can robustly reconstruct complex dynamics (Part 1, Fig. \ref{fig:flow-chart}). After fitting each window, we recover a collection of DMD models describing a range spatiotemporal dynamics (Part 2, Fig. \ref{fig:flow-chart}). The temporal dynamics of the DMD models are clustered into bands. The high frequency components are well-resolved and removed. The low-frequency component is not well-resolved and is withheld as the input for the next decomposition level with a larger window, which can now resolve more of the low-frequency component of the signal. Finally, as information tends to leak between levels, a final global determination of the temporal bands is made, yielding an interpretable decomposition capable of resolving non-stationary, transient, multi-scale, and high-dimensional features (Part 3, Fig. \ref{fig:flow-chart}).

To demonstrate the ability of mrCOSTS, we diagnose three systems characterized by complex multi-scale dynamics. We choose to demonstrate mrCOSTS on real systems instead of toy models as the data from toy models do not approximate the true complexities of multi-scale data. That said, an example with a toy model can be found as part of the PyDMD tutorial of mrCOSTS \cite{ichinagaPyDMDPythonPackage2024}. In the first example, we diagnose patterns of sea surface temperature (SST) over the central Pacific over a 150 year period, focusing on the 2015-2016 extreme SST anomaly. In the second example, we examine Local Field Potential (LFP) observations of the motor cortex of a monkey during a trained grasping task. Finally, we decompose ground-based remote sensing retrievals of horizontal wind speed observing the confluence of a tributary and main valley flow in the mountain boundary layer (MoBL), which is well-known for its complex multi-scale dynamics.

%%%%%%%%%%%%%%
%%%%%%%%%%%%%%
%%%%%%%%%%%%%%
\section*{Multi-resolution Coherent Spatio-Temporal Scale Separation (mrCOSTS)}\label{sect:mrCOSTS-description}
mrCOSTS builds on the principles of multi-resolution Dynamic Mode Decomposition \cite[mrDMD;][]{Kutz2016_dmd}, specifically the sliding mrDMD method \cite{Dylewsky2019_dmd}, but is tailored to the complexities of multi-scale data.

\begin{figure*}[t!]
\centering
\includegraphics[width=17.8cm]{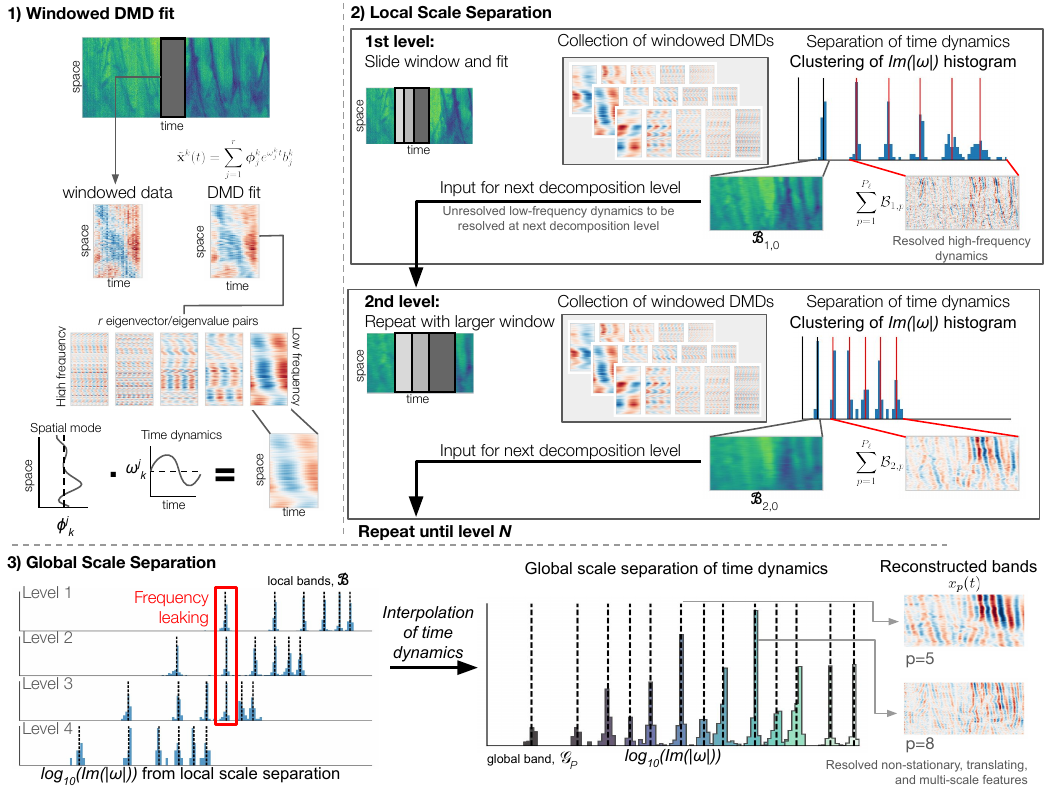}
\caption{An overview of the mrCOSTS algorithm. (1) The data are windowed and a DMD model is fit to each. The resulting fit is comprised of several time dynamics with distinct frequencies with a corresponding spatial mode, which can take on any functional shape. For the local scale separation (2) the window is slid across the data resulting in a collection of DMD fits. The time dynamics are clustered yielding bands comprised of eigenvalue/eigenvector pairs which share time dynamics. The bands with the highest frequency are removed and the data reconstruction using the lowest frequency band is given to the next decomposition level. At the next level the window size is increased, allowing slower time dynamics to be fit, and the process is repeated. (3) Finally, the collection of frequency components from the local scale separation need to be merged to account for frequency leaking prior to a final, global scale separation. The resulting bands decompose the complex multi-scale physics of the input signal.}
\label{fig:flow-chart}
\end{figure*}%[tbhp]

DMD is a method for data-driven model discovery \cite{Rowley2009_dmd, Kutz_DMD-textoobk_2016, schmidDynamicModeDecomposition2022} that provides a low-rank, interpretable model of coherent spatiotemporal features. However, a basic assumption in many variants of the DMD \cite[for a list of commonly applied variants see][]{ichinagaPyDMDPythonPackage2024} is that the data can be approximated by a (relatively) small number of coherent spatiotemporal modes. For real world systems, especially those characterized by multi-scale dynamics, this assumption on the time dynamics is not often valid. Consequently, DMD analysis relies on strategies such as substantial window crafting \cite{chenVariantsDynamicMode2012, bagheriKoopmanmodeDecompositionCylinder2013} or using steady-state systems \cite{Taira2020}. In response, the multi-resolution DMD (mrDMD) was developed \cite{Kutz2016_dmd, Dylewsky2019_dmd}, however, we found these methods, while sufficient for synthetic test data, struggled with real multi-scale data leading to the development of mrCOSTS.

The mrDMD framework is given by Eq. \ref{eq:windowed-dmd}, in which input snapshots $\mathbf{x}(t) \in \mathbb{R}^n$ collected across times $\{t_i\}_{i=1}^m$ are sequentially windowed and fit with a DMD model so that
\begin{equation}
\tilde{\mathbf{x}}^k(t)=\sum^{r}_{j=1}\boldsymbol{\phi}^k_j e^{\omega^k_jt}b^k_j + \mathbf{c}_k\label{eq:windowed-dmd}. \\
\end{equation}
We use $\tilde{\mathbf{x}}^k(t)$ to indicate an approximation to the windowed data, $\mathbf{x}^k(t)$. We use $k$ to index the data windows so that snapshots belonging to the $k$th window, denoted as $\mathbf{x}^k(t)$, are approximated by the decomposition given by Eq. \ref{eq:windowed-dmd}.

The input data can have arbitrary dimensionality but the non-time dimensions must be flattened to a single spatial dimension, similar to the data shape requirement for most machine learning algorithms \cite{scikit-learn2011}. Similarly, $\mathbf{x}(t)$ may be composed of multiple variables, as demonstrated in the mountain boundary layer example later. In the case of multiple variables, the DMD model finds the coherent spatiotemporal patterns shared among the variables.  

Prior to decomposition, each window is normalized by removing the time mean, denoted $\mathbf{c}_k$, which must be added back in to reconstruct $\mathbf{x}^k(t)$ \citep{Hirsh2020}. We use $j$ to index the DMD eigenvalue ($\omega$)/eigenvector ($\boldsymbol{\phi}$) pairs up to rank $r$, where $b$ gives the amplitude of the pair. The spatial structure of the data, which is encoded by $\boldsymbol{\phi}$, has time dynamics that are determined by $\omega$, where the real component of $\omega$ governs exponential growth or decay and the imaginary component governs oscillatory behaviors (part 1 Fig. \ref{fig:flow-chart}). Eq. \ref{eq:windowed-dmd} is then fit to sliding windows of the data with a fixed window length in time, which constitutes a decomposition level. Thus for a given decomposition level, we obtain a collection of DMD fits that describe the coherent spatiotemporal structures of the data (part 2 Fig. \ref{fig:flow-chart}).

Since the decomposition given by Eq. \ref{eq:windowed-dmd} differs across all windows $k$ for a given decomposition level, indexed with $\ell$, we cluster all $Im(|\omega^k_j|)$ values from decomposition level $\ell$. This yields clustered eigenvalue/eigenvector bands, denoted $\mathcal{B}$, that are local to the $\ell$th decomposition level (part 2 Fig. \ref{fig:flow-chart}). $\mathcal{B}_{\ell, 0}$ has the smallest $Im(|\omega|)$, i.e. the slowest frequencies, and is characterized by time scales that are longer than the window length. This band is usually decomposed into more refined time scales at larger decomposition levels and consequently $\mathcal{B}_{\ell, 0}$ is not well-resolved. In contrast, the bands with larger $Im(|\omega|)$, i.e. higher frequencies, which we denote in order of increasing frequency using $\mathcal{B}_{\ell, 1}, \mathcal{B}_{\ell, 2}, \dots, \mathcal{B}_{\ell, P_\ell}$, are better resolved at the $\ell$th decomposition level. The contribution of any band $\mathcal{B}_{\ell, p}$ to the signal can be recovered by
\begin{equation}
    \tilde{\mathbf{x}}_{\ell, p}(t) = \sum_{k} \sum_{j \in \mathcal{B}_{\ell, p}} \boldsymbol{\phi}^k_j e^{\omega^k_jt}b^k_j + \mathbf{c}_k.\label{eq:local-recon}
\end{equation}
We thus use Eq. \ref{eq:local-recon} for $p=0$ to reconstruct the unresolved low-frequency components of the data at decomposition level $\ell$ so that it may be used as input for the next decomposition level. The background value, $\mathbf{c}_k$, is only included for $p=0$ and is otherwise excluded.

The algorithm is then repeated for decomposition level $\ell + 1$, which uses a larger data window with a length selected by the user (see Hyperparameters). In doing so, high-frequency components which are resolvable by DMD at decomposition level $\ell$ are removed, while the low-frequency components are retained for further scale separation at larger window sizes. The algorithm is repeated until the largest desired decomposition level, $N$, is reached.

The general algorithm outlined above, which follows the broad outline of the mrDMD \cite{Dylewsky2019_dmd}, is called the local scale separation. A normalizing factor was omitted for reconstructing the overlapping windows for brevity. The windowed DMD fit is most robust near the center of the window and less reliable near the edges \citep{Dylewsky2019_dmd}. Overlapping DMD windows are strongly weighted towards the center of the window as a consequence of the poor fitting at window edges. For this reason, mrCOSTS is least reliable at the edges of the time domain. This error is analogous to the cone-of-influence from continuous wavelet transformations \citep{Lau1995} and is thus most noticeable at the largest decomposition level with the longest time scales. 

Performing this algorithm recursively on data without clearly separable time dynamics, i.e., multi-scale data, can lead to poor results due to the leaking of information across decomposition levels (part 3 in Fig. \ref{fig:flow-chart}) as well as other deficiencies in the original algorithm (see Methods). To remedy the leaking of frequency components between levels, a global scale separation has to be performed on the fast frequencies from all decomposition levels. For this step, the $\omega$ components from all of the local frequency bands $\mathcal{B}_{\ell, p}$ for all $\ell$ and $p > 0$ are interpolated to the mean time of $\mathcal{B}_{0, p}$'s windows using a nearest neighbor approach. The interpolated $\omega$ are then used to perform the global scale separation and are then re-indexed back to the original decomposition window's time windows. The resulting globally scale separated bands, $\mathcal{G}_{0}, \mathcal{G}_{1}, \dots, \mathcal{G}_{P}$, listed in order of increasing band frequency, are a complete scale separation that includes the leaked component between each level. We again use $p$ to index the globally scale separated bands.

We can then reconstruct the contribution of $\mathcal{G}_{p}$ to $\tilde{\mathbf{x}}_p(t)$ using
\begin{equation}
\tilde{\mathbf{x}}_p(t)=\sum_{k}\sum_{(j, \ell) \in \mathcal{G}_p} \boldsymbol{\phi}^k_{j,\ell} e^{\omega^k_{j,\ell}t}b^k_{j,\ell} \label{eq:band-reconstruction}. \\
\end{equation}
A global reconstruction of the mrCOSTS fit can be found by 
\begin{equation}
\tilde{\mathbf{x}}(t)\approx\sum^{}_{p} \tilde{\mathbf{x}}_p(t) +  \mathbf{c}_{P, k} \label{eq:global-reconstruction}. \\
\end{equation}
The background value of each window in $\mathcal{G}_{0}$ (the lowest frequency band), $\mathbf{c}_{0, k}$, needs to be included for the global reconstruction of the signal. Note, $\mathcal{G}_0$ is defined as the low-frequency component of the largest decomposition level and therefor may not contain well-resolved dynamics. The influence of multiple selected bands can be found through Eq. \ref{eq:global-reconstruction} by simply summing over the desired $p$ and omitting the background value. For ease of interpretability we identify $\mathcal{G}$ bands according to the cluster centroid period or frequency (vertical dashed lines in part 3 Fig. \ref{fig:flow-chart}), depending on which is more intuitive for the system of interest. 

A useful property of mrCOSTS is that each $\mathcal{G}_p$ for $p>$0 is an approximation of the contribution by coherent spatiotemporal processes with a narrow frequency range to the fluctuating component of Reynolds decomposition. This property makes mrCOSTS especially powerful for the decomposition of systems governed by the Navier-Stokes equations as well as other systems with similar multi-scale hierarchies. Another important property is that non-spatially coherent processes, such as white noise, cannot be fit by mrCOSTS as white noise is not spatially coherent with an oscillatory time scale. As a result, $\tilde{\mathbf{x}}$ is an approximation of the de-noised $\mathbf{x}$.

%%%%%%%%%%%%%%
%%%%%%%%%%%%%%
%%%%%%%%%%%%%%
\section*{Results}

As stated previously, we demonstrate mrCOSTS on three examples. Each example is characteristic of the multi-scale dynamics which generally frustrate other analysis methods and are presented in order of our subjective assessment of complexity as well as unknown dynamics. The SST example provides a validation of sorts due to the well-known time scales governing the system. Using mrCOSTS we recover these known time scales while also recovering the generally poorly-described spatial patterns as well as previously unknown time scales critical to the description of this system. In the neurology example, we trivially isolate and characterize well-known frequency bands and in doing so isolate the previously poorly-described translating spatial patterns. Finally, in the MoBL example we recover entirely unknown dynamics for a system that generally frustrated rigorous and objective analysis.

\subsection*{Sea surface temperature}

\begin{figure*}[t!]
\centering
\includegraphics[width=17.8cm]{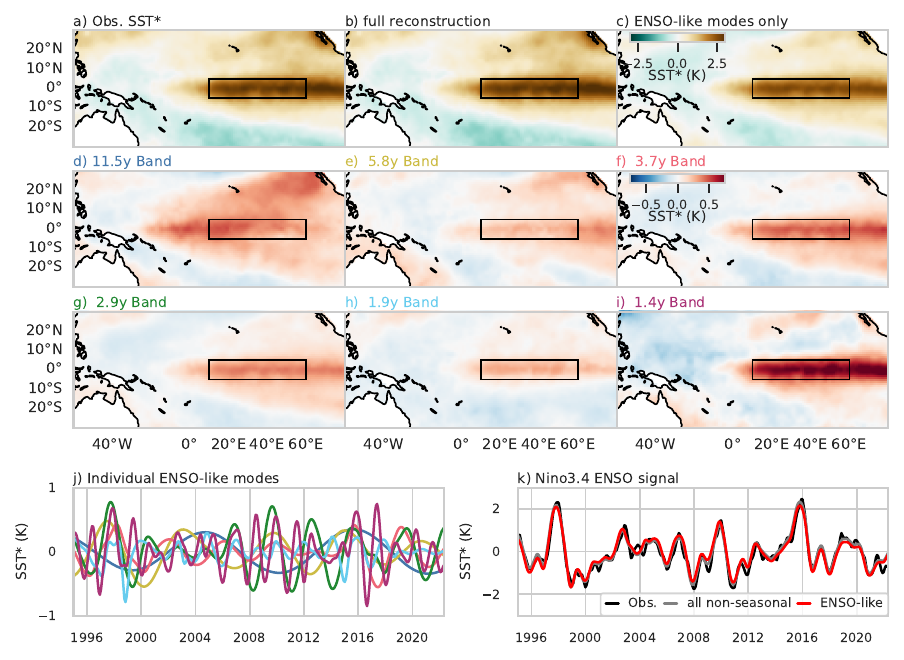}
\caption{The mrCOSTS diagnosis of the ENSO anomaly during the period Dec 2015 - Feb 2016. These dynamics were recovered with no hyperparameter tuning. Shown are the (a) observed SST$^*$, (b) mrCOSTS reconstruction of SST$^*$ using all modes, and (c) mrCOSTS reconstruction using only the ENSO-like modes. (d-i) Each of the ENSO-like modes is reconstructed individually with the band centroid indicated in the subplot label. The Ni\~no 3.4 box is drawn for reference. (j) The spatial mean of each of the ENSO-like modes inside the Ni\~no 3.4 box are shown, yielding a time series of each mode's contribution to SST$^*$. Each mode is color coordinated with the subplot titles in (d-i) . (k) The time series of the observed, ENSO-like modes, and all non-seasonal modes SST* within Ni\~no 3.4 are compared. The complex time evolution of ENSO is recovered and explained in no small part due to discovering the previously undetected time dynamics of the 11 y, 1.9 y, and 1.4 y bands.}
\label{fig:ENSO-reconstruction}
\end{figure*}%[tbhp]

Climatic SST patterns of the Pacific Ocean, especially the El Ni$\tilde{n}$o-Southern Oscillation, are major modes of variability in the global atmospheric circulation, driving processes in the earth system with substantial societal impacts \cite{glantzReviewingOceanicNino2020, callahanPersistentEffectNino2023}. Predicting and understanding the dynamics of ENSO is a major research challenge \cite{mcphadenENSOIntegratingConcept2006, santosoDefiningCharacteristicsENSO2017}, especially in the context of possible changes to this mode of internal variability with climate change \cite{caiChangingNinoSouthern2021, callahanPersistentEffectNino2023}.

It is well known that quasi-periodic behavior of the multi-scale processes present in SST data are a substantial challenge for analysis \cite{dommengetCautionaryNoteInterpretation2002, dommengetEvaluatingEOFModes2007}. ENSO is thought to have a dominant period between 2 to 7 years \cite{Torrence1998, Torrence1999} with a chaotic or complex quasiperiodic behavior \cite{tzipermanIrregularityLockingSeasonal1995}.

A wide variety of metrics for ENSO have to be used for various practical and research reasons \citep{bamstonDocumentationHighlyENSO1997, trenberthDefinitionNino1997, trenberthIndicesNinoEvolution2001, lheureuxObservingPredicting20152017, santosoDefiningCharacteristicsENSO2017, adamsonSituatingNinoCritical2022, trenberthkevinNinoSSTIndices, noaacliamtepredictioncenterNOAAClimatePrediction}. To describe the mrCOSTS decomposition of SST, we use the Ni\~no 3.4 box \cite{bamstonDocumentationHighlyENSO1997, trenberthDefinitionNino1997}, which is a result of careful window crafting, but perform the decomposition on a much larger region. ENSO is commonly described using monthly SST anomalies, which we denote as SST$^*$ (see Methods).

To demonstrate the unsupervised diagnostic power of mrCOSTS, no hyperparameter tuning was performed and no window crating was applied. The mrCOSTS fit recovers the complex time dynamics of the SST observations, especially those at time scales longer than the annual cycle including the time scales commonly associated with ESNO (Fig. S1). The resulting $\mathcal{G}_p$ were visually inspected to identify a set of 6 bands with ENSO-like spatial patterns. The time scales of these modes varied from 1.4 years to 11 years. Using these bands, the extreme ENSO anomaly in the period of DJF 2015-2016 \citep{lheureuxObservingPredicting20152017, santosoDefiningCharacteristicsENSO2017} was reconstructed using Eq. \ref{eq:global-reconstruction} (Fig. \ref{fig:ENSO-reconstruction}b) and Eq. \ref{eq:band-reconstruction} for the ENSO-like bands alone (Fig. \ref{fig:ENSO-reconstruction}c). The mrCOSTS decomposition accurately reconstructs the domain-wide SST$^*$ (Fig. \ref{fig:ENSO-reconstruction}a,b) with most of the variability within the domain being reconstructed by the ENSO-like bands (Fig. \ref{fig:ENSO-reconstruction}a-c). The ENSO-like bands almost entirely decompose the SST$^*$ signal in the Ni\~no 3.4 box with the reconstruction using all supra-seasonal bands barely altering the reconstruction (Fig. \ref{fig:ENSO-reconstruction}k).

The six ENSO-like bands include 3 bands outside the time scales normally associated with ENSO: 11 years, 1.9 years, and 1.4 years. The 2015 anomaly was the result of all 6 ENSO-like bands having a positive anomaly at the same time (Fig. \ref{fig:ENSO-reconstruction}j), a remarkably rare occurrence in the reconstruction, not occurring for any other positive ENSO event. The largest contributor to the 2015 anomaly was the 1.4 year band, a time scale neglected in the characterization of ENSO, potentially connected to the known but difficult to characterize relationship between extreme ENSO events and the seasonal cycle \citep{santosoDefiningCharacteristicsENSO2017}. Other strong positive ENSO anomalies only had contributions from 4 or 5 of the ENSO-like bands, partially explaining either their weaker expression in the Ni\~no 3.4 box or differences in spatial patterns and temporal evolution relative to the 2015-2016 event \citep{trenberthIndicesNinoEvolution2001, lheureuxObservingPredicting20152017, santosoDefiningCharacteristicsENSO2017}. For instance, the positive SST$^*$ in 1997 was a result of all ENSO modes in the positive phase except for the 11.5 year band (Fig. \ref{fig:ENSO-reconstruction}j). Further, the distinct ENSO-like modes are remarkably well-separated given that mrCOSTS was fit over the entire central Pacific Ocean over a 150 year long period. Previous analyses could not use the entire data set to trivially recover these results, but would need to do either some spatial aggregation \cite[e.g,][]{Torrence1998} or window crafting \cite[e.g.,][]{bamstonDocumentationHighlyENSO1997}.

Unsurprisingly, given the known complexity and quasiperiodic nature of the ENSO signal, the individual ENSO-like bands go through periods of varying activity and do not exhibit clean oscillatory patterns. The 11 year band turns on in approximately 1980 but becomes surprisingly regular into present day afterwards (dark blue line, Fig. \ref{fig:ENSO-reconstruction}j). In contrast, the 1.9 and 1.4 year bands exhibit little regularity (light blue and purple lines, Fig. \ref{fig:ENSO-reconstruction}j). Despite these non-stationary processes mrCOSTS was able to characterize and trivially diagnose the chaotic behavior of ENSO, simultaneously recovering the time dynamics and their associated coherent spatial patterns. While we recovered the dynamics known to exist in this system, we additionally identified three frequency bands which were otherwise not recognized for describing the system. Using these results we highlight how mrCOSTS was able to diagnosis in an unsupversived fashion known complex multi-scale dynamics while also explaining dynamics that were previously difficult to characterize.

\subsection*{Neurology}

\begin{figure*}[t!]
\centering
\includegraphics[width=17.8cm]{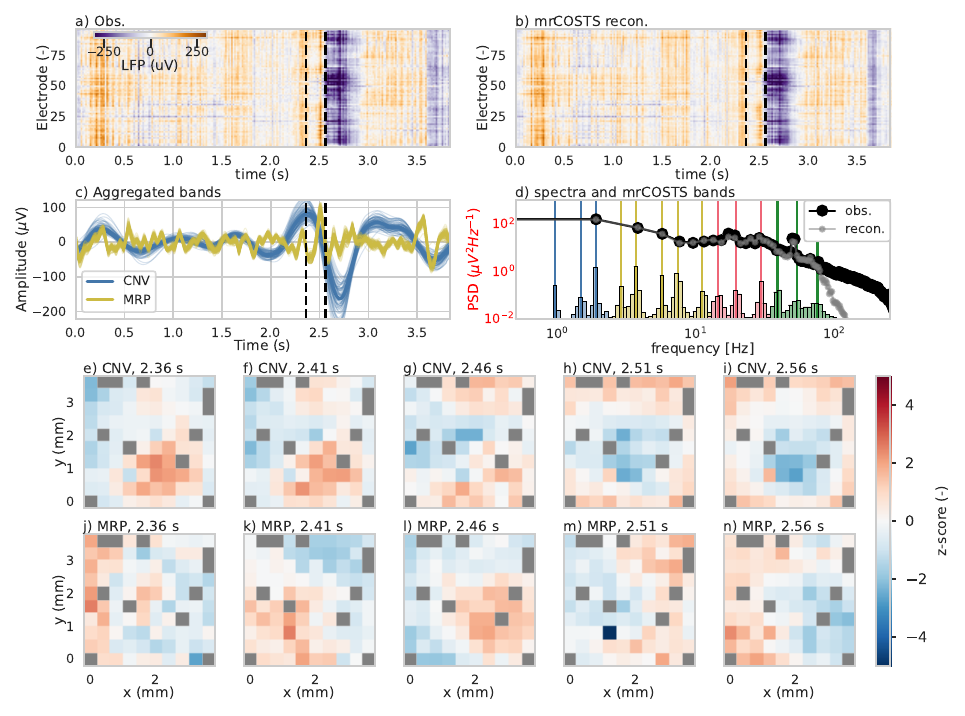}
\caption{Decomposition of the LFP signal to the CNV and MRP bands using mrCOSTS for a selection of times. See Movie S1 for all bands and times. (a) Observed LFP and the (b) mrCOSTS reconstruction. (c) The aggregated bands were constructed to be consistent with the CNV (blue) and MRP (yellow) bands using Eq. \ref{eq:band-reconstruction}. These bands also shared characteristics during the movement. Each electrode from (a) and (b) is plotted as a separate line in (c). In (a-c) the first vertical dashed line denotes the cue for movement ($t_1$=2.36 s) while the second line indicates the beginning of the movement ($t_2$ = 2.51 s). (d) The mean power spectral density of the LFP (black line) and the mrCOSTS reconstruction (grey line). Plotted underneath is the histogram of the globally interpolated $\omega$ on a relative scale. The histogram is colored coordinated with the assigned bands shown in (c) for the CNV and MRP and SI movie 1 for all bands. The cluster centroids are indicated as vertical lines. The CNV (e-i) and MRP (j-n) bands are spatially reconstructed for five times from $t_1$ to $t_2$ with an even time spacing. The spatial reconstructions are standardized by a z-score to facilitate their plotting on the same color scale. Grey points indicate no data. The translating spatial patterns of these bands demonstrates how mrCOSTS can decompose these types of complex spatiotemporal signals, which are otherwise not amenable to analysis.}
\label{fig:monkey}
\end{figure*}%[tbhp]

Electrophysiological observations of the brain are thought to be a method for uncovering the basis of cognition in neural processes \citep{einevollModellingAnalysisLocal2013}. The high-frequency component of the signal ($>$500 Hz) is related to the spiking of individual neurons near the electrode while the low-frequency component of the signal ($<$ 500 Hz; the local field potential, or LFP) offers an insight into the integrated behavior of the entire neural circuit, albeit with substantial complications to this general interpretation \citep{herrerasLocalFieldPotentials2016a}.

The relationship between LFP signals and specific neural pathways is complex, as the contributions to the LFP signal from specific neural circuits are ambiguous and are comprised of complex spatiotemporal signals \citep{einevollModellingAnalysisLocal2013, herrerasLocalFieldPotentials2016a}. Temporal dynamics occur across a range of scales and are activated by different spatial regions/neural circuits and sometimes even the same set of neurons. Consequently, the ``spatial factors [of LFPs] are among the most important issues to be explored in depth in the following years" \citep{herrerasLocalFieldPotentials2016a}. Developing mathematical tools for analyzing the LFP is hence seen as one of the major needs in neuroscience \citep{einevollModellingAnalysisLocal2013, riehleMappingSpatiotemporalStructure2013}. DMD has been found as a promising approach for analyzing LFP data \citep{bruntonExtractingSpatialTemporal2016, ferreNonStationaryDynamicMode2023}, but the multi-scale nature of these data frustrates the application of DMD more generally. From this perspective, we therefore reframe the problem of interpreting LFPs as a challenge, in large part, in interpreting data with multi-scale dynamics.

The LFP signal is often described in terms of frequency bands, making the scale-separation task particularly amenable to decomposition with mrCOSTS. The observations we use cover frequencies describing contingent negative variation (CNV) related to delayed tasks \citep[$<$ 2 Hz;][]{zaepffelPlanningVisuallyGuided2012}, the movement-related potential (MRP) in the motor cortical areas which encodes movement type and direction \citep[3-15 Hz;][]{rickertEncodingMovementDirection2005, kilavikEvokedPotentialsMotor2010}, beta waves \citep[15-30 Hz;][]{niedermeyerElectroencephalographyBasicPrinciples2005, kilavikContextRelatedFrequencyModulations2012}, and gamma waves\citep[$>$30 Hz;][]{niedermeyerElectroencephalographyBasicPrinciples2005}. The MRP is of particular interest due its relationship with specific movements through complex spatiotemporal patterns \citep{mehringInferenceHandMovements2003, rickertEncodingMovementDirection2005, markowitzOptimizingDecodingMovement2011, riehleMappingSpatiotemporalStructure2013}. This relationship is particularly difficult to describe, yet highly desirable to understand for applications such as neuroprostheses \citep[e.g.,][]{markowitzOptimizingDecodingMovement2011}. Often the MRP is interpreted through intuitive transformations such as the aggregatig a single electrode across multiple trials or aggregating the spatial components, neglecting aspects of the multi-scale dynamics.

We use electrophysiological observations of a macaque's motor cortex during a trained exercise \citep{riehleMappingSpatiotemporalStructure2013}. A monkey, with chronically implanted 10-by-10 Utah electrode array, was trained for a delayed reach-to-grasp task consisting of a cue for the type of grip (side or precision), a wait interval, and a go signal encoding how hard to grip (low- or high-force). Upon successful completion a reward was dispensed. An advantage of these data is that they were released for the explicit purpose of developing analytical techniques for electrophysiological observations \citep{brochierMassivelyParallelRecordings2018}.

We applied mrCOSTS to trial two from monkey L (a side-grip low-force trial). The original hyperparameters successfully decomposed the system, but we found that a very low-frequency oscillation ($\approx$1 s) was missed by the initial set of hyperparameters. Consequently, the hyperparameters were adjusted exactly once to include an additional decomposition level with a very large window size explicitly to capture this time scale. These very-low frequency features are difficult to characterize \textit{a priori} using typical methods like the power spectral density due to their poor resolution for very long time scales (Fig. \ref{fig:monkey}d). However, mrCOSTS easily identifies the slower time dynamics of the system, demonstrating some of the diagnostic power of mrCOSTS over traditional techniques (compare the histogram of $\omega$ and band centroids to the power spectral density in Fig. \ref{fig:monkey}d).

The reconstruction using the mrCOSTS fit captures all of the salient features of the observations (Fig. \ref{fig:monkey}a, b, and d). A total of 14 bands were identified in the global scale separation (Fig. \ref{fig:monkey}d). We aggregated the mrCOSTS $\mathcal{G}_p$ into four bands using Eq. \ref{eq:band-reconstruction} to align with the CNV, MRP, beta oscillations, and gamma oscillations. All bands and their spatial structure are shown in Movie S1, whereas only the MRP and CNV are presented in Fig. \ref{fig:monkey} for visual clarity. The mrCOSTS fit recovers the spectra for the time scales slower than $\approx$ 75 Hz, meaning the dynamics of all of the neurological bands, especially the MRP and CNV, are well resolved (Fig. \ref{fig:monkey}d).

With mrCOSTS, we trivially retrieve the spatial patterns of the MRP and CNV individually, which have distinct behavior during the grasping movement and reveal complex traveling spatiotemporal patterns (Movie S1, Fig. \ref{fig:monkey}e-n). The CNV aggregated band (blue line Fig. \ref{fig:monkey}c) has a strong positive anomaly in the lower right hand corner of the electrode array (Fig. \ref{fig:monkey}e) which travels out of frame to the bottom right corner and is replaced with a strong negative anomaly in roughly the same location (Fig. \ref{fig:monkey}e-i, Movie S1). The MRP aggregated band also features translating processes. The positive anomaly in the lower left hand side of the array at 2.4 s travels to the right across the array (Fig. \ref{fig:monkey}j-n, Movie S1). The MRP and CNV signals are not spatially static oscillations and appear to characterized by these types of spatially translating features (Movie S1).

These sort of translating patterns of the MRP and CNV and their potential relationship with specific movements or events, respectively, are not well-explored even though the spatial patterns are known to be important \citep{herrerasLocalFieldPotentials2016a}. We speculate these patterns were not explored as a consequence of the multi-scale nature of the data. The cursory analysis demonstrates some of the potential of unsupervised and structured diagnoses of mrCOSTS for multi-scale data. 

\subsection*{Mountain Boundary Layer (MoBL)}

\begin{figure*}[t!]
\centering
\includegraphics[width=17.8cm]{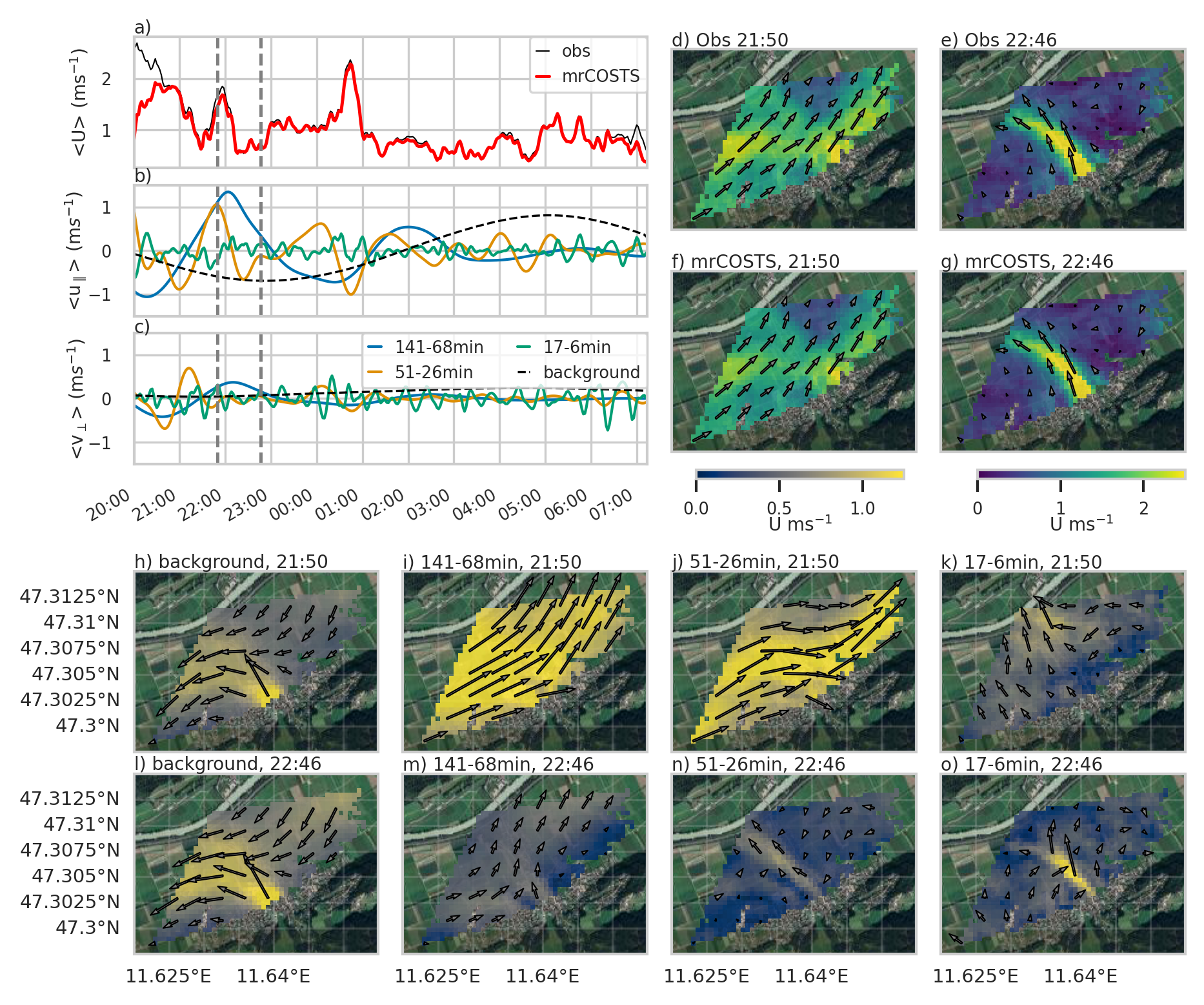}
\caption{An overview of the mrCOSTS decomposition of the MoBL data highlighting two times of interest, see Movie S2 for all times. These dynamics were found with little hyperparameter tuning in a completely unsupervised fashion. (a) A comparison of the observed and mrCOSTS reconstruction of the mean wind speed. The time series of the aggregated bands for (b) the along-valley wind $\langle u_{\parallel} \rangle$, (c) the across-valley wind $\langle v_{\perp} \rangle$ are shown. (d,e) The co-planar retrieval of wind speed and direction approximately 60 m above the valley floor at two times with (f,g) the mrCOSTS reconstructions. The two times are indicated as vertical dashed lines in (a-c). The reconstructed mean wind speed and direction for the background mode and the three aggregated bands are shown for the first (h-k) and second (l-o) times of interest. The mrCOSTS decomposition reveals the complex interaction of the background, seich-like oscillations, and tributary valley in-flow which could not be well-characterized previously. The background image \textcopyright 2024 CNES/Airbus, Google, Maxar Technologies.}
\label{fig:MoBL-reconstruction}
\end{figure*}%[tbhp]

The MoBL is well-known for its multi-scale dynamics \citep{Rotach2017, Serafin2018, rotachCollaborativeEffortBetter2022}, with features such as frustratingly difficult to analyze non-monochromatic waves \citep{nappoIntroductionAtmosphericGravity2013, Sun2015}, overlapping process scales \citep{Lehner2018, pfisterTEAMxPC22Alpine2024}, and dynamics occurring across orders of magnitude \citep{zardiDiurnalMountainWind2013, Lehner2018}, all with the additional constraint of limited observations \citep{rotachCollaborativeEffortBetter2022, pfisterTEAMxPC22Alpine2024}. These difficult to characterize multi-scale dynamics form a fundamental obstacle to the modeling, understanding, and prediction of the atmosphere over regions of complex terrain.

The MoBL is an illustrative example of the challenges of discovering unknown dynamics from multi-scale data. Atmospheric models are developed by defining a grid spacing. The size of the grid spacing dictates which processes are explicitly resolved (those larger than the grid spacing) and those which must be parameterized \cite[those smaller than the grid spacing;][]{leonardEnergyCascadeLargeEddy1975, wyngaardSurfaceLayerSurface1998, honnertAtmosphericBoundaryLayer2020}. The construction of these parameterizations relies on understanding how processes at a smaller scale create emergent behavior at a larger scale, which has to be informed by data \cite{lemone100YearsProgress2019, edwardsRepresentationBoundaryLayerProcesses2020a}. As the grid spacing shrinks, i.e. when seeking to improve the model, new parameterizations have to be derived \cite{chowCrossingMultipleGray2019a, honnertAtmosphericBoundaryLayer2020}. These improvements require richer, more detailed data, which due to their multi-scale nature come up against the limits of traditional analytical techniques \cite{pfisterTEAMxPC22Alpine2024}. Thus, even though the physical laws of fluid flows are well-known, describing the emergent physics across scales remains a fundamental need and challenge, in part because of the multi-scale nature of MoBL data.

An observational campaign was conducted in the summer of 2022 in the Inn Valley, Austria for the explicit purpose of addressing these challenges \citep{pfisterTEAMxPC22Alpine2024}. We make use of the co-planar retrieval of horizontal wind speed components, $u$ (wind along the east-west axis) and $v$ (wind along the north-south axis), using two horizontally-oriented LIDARs with overlapping fields of view \citep{Kalthoff2014}. The co-planar observations were set up to observe the outflow of the tributary Weer Valley into the main Inn Valley 60 m above the valley floor. These types of data are generally not amenable to objective analytical methods and instead require careful human interpretation. The individual wind speed components can be used to recover the mean horizontal wind speed $U=\sqrt{u^2 + v^2}$ and direction. We use $<>$ to denote the spatial mean. The wind speed components were rotated into a coordinate system better aligned with the Inn Valley axis (roughly southeast to northwest), yielding an across-valley component, $v_{\perp}$ and an along-valley component, $u_{\parallel}$.

The confluence of the two valley flows are inherently a difficult multi-scale challenge. The tributary flow is characterized by scales an order of magnitude smaller than the main valley flow, with the scales of tributary flow likely overlapping with the process scales of the main valley. We say ``likely" as this MoBL process eluded prior systematic and objective analysis \citep{zardiDiurnalMountainWind2013} and one must instead hypothesize the process scales based on the valley geometry \citep{pfisterTEAMxPC22Alpine2024}. Additionally, the MoBL processes are known to be non-stationary and translate, subject to largely unknown dynamics and forcings, with observations corrupted by noise. For this reason, we subjectively characterize the MoBL data as the most difficult example of multi-scale dynamics examined in this study.

As in the neurology example, we performed a small amount of hyperparameter tuning in order to capture some very-low frequency dynamics which were only apparent after performing an initial mrCOSTS decomposition. This was achieved simply by performing an additional decomposition with a much larger window size.

From the individual $\mathcal{G}_p$ we formed 3 aggregate bands of 141 min to 68 min, 51 min to 26 min, and 17 min to 6 min using Eq. \ref{eq:band-reconstruction}. We also highlight the background mode that evolves very slowly on an $\approx$ 10 hour time scale, which was too long to decompose given the duration of the data (9 hours) but whose temporal dynamics are fairly clean, with a period that appears uncontaminated by other time scales (Fig. \ref{fig:MoBL-reconstruction}b). 

The mrCOSTS decomposition reveals a remarkable oscillation of the wind along the valley axis. These oscillations were most prominent in 95 min and 51 min bands and these two bands dominate the aggregated band behavior (blue line and orange line in Fig \ref{fig:MoBL-reconstruction}b, respectively). The clean and regular oscillations of the along-valley wind component, especially in the non-aggregated bands (not shown), are suggestive of a standing wave, as regular external forcing of the valley system on these time scales is less likely. A potential candidate process is a standing wave known as a seiche, a process analogous to water sloshing in a bath tub. These phenomenon are most well-known in enclosed basins in the atmosphere \citep{whitemanMetcrax2006Meteorological2008, zardiDiurnalMountainWind2013}, where the geometry of the basin is regular, potentially allowing for clean harmonics of the seiche to be more easily determined both analytically and visually. However, the scale of these waves is unclear, as their spatial scales are larger than the observation footprint. Further, the boundary upon which the waves would reflect is also unclear, but potentially could be related to the bend in the valley system immediately to the east of the observations or from persistent density boundaries caused by the multiple tributary in-flows along the main valley.

The presence of these apparently standing waves, which we will refer to as seiche-like for simplicity, in a large valley system like the Inn Valley is a novel result not previously recognized. We speculate this is because of the imperfect, multi-scale oscillations interacting with the other dynamics of the MoBL, which include a wide-range of processes acting in three dimensions. Extracting out this pattern was only possible using the mrCOSTS decomposition. Phenomenon like gravity waves are not explored on these horizontal scales generally, in no small part due to a lack of analytical tools for non-monochromatic waves \citep{nappoIntroductionAtmosphericGravity2013}, a constraint easily overcome with mrCOSTS.

We highlight two times of interest to  demonstrate how the mrCOSTS decomposition reveals the complex dynamics of this system ($t_1$=2150 and $t_2$=2246; vertical dashed lines in Fig. \ref{fig:MoBL-reconstruction}a-c) while Movie S2 shows the entire time period. At $t_1$ the tributary outflow appears to not exist with a substantial down valley flow in its place (Fig. \ref{fig:MoBL-reconstruction}d,f). However, with mrCOSTS we can reveal that the very slowly fluctuating background containing a persistent tributary outflow and main valley flow is overwhelmed by the seiche-like motions (Fig. \ref{fig:MoBL-reconstruction}h vs \ref{fig:MoBL-reconstruction}i,j). The highest frequencies reveal a cross valley oscillation (Fig. \ref{fig:MoBL-reconstruction}c,k). In other words, the dynamics of the tributary valley in-flow was relatively constant in time (Fig. \ref{fig:MoBL-reconstruction}h) but the seiche-like oscillations along the valley floor masked this behavior (Fig. \ref{fig:MoBL-reconstruction}i,j). 

For $t_2$, the slowly fluctuating background mode is largely unchanged from $t_1$ except for a strengthening of the up-valley wind (Fig. \ref{fig:MoBL-reconstruction}l). The seiche-like oscillations of the largest time scales (Fig. \ref{fig:MoBL-reconstruction}m) oppose the main valley flow, while the finer time scales (Fig. \ref{fig:MoBL-reconstruction}n,o) reveal a pulse in the tributary flow, which travels into the main valley over the next few time steps (Movie S2). The composite effect is to create the appearance of a strong tributary flow penetrating into a weak main valley. The highest frequencies characterize the sharp boundary between the valley flow and the tributary flow injecting into it, however, the pulsing of the tributary flow is rare and does not characterize the flow generally (Movie S2). Instead, the interplay between the nearly constant background flow and the seiche-like motions dominates the system, even in periods in which we would have interpreted the tributary in-flow as acting alone such as during $t_2$. Additionally, the background up-valley flow (negative $u_{\parallel}$ in Fig. \ref{fig:MoBL-reconstruction}b) is unusual and counter to the prevailing understanding of thermally driven flows in complex terrain \citep{zardiDiurnalMountainWind2013}. 

Flow patterns such as those highlighted here could not be previously extracted due to the lack of tools for analyzing these kinds of multi-scale data. Being able to rigorously diagnose these kinds of multi-scale data will enable us to refine and improve our intuitive interpretation of these data as well as develop new physical explanations. Revealing these kinds of hidden dynamics in multi-scale data further underlines the utility offered by an unsupervised diagnostic like mrCOSTS.

\section{Discussion}\label{Sec:Discussion}

In summary, we have provided a new definition for multi-scale data for which existing computational methods are inadequate. To answer the challenge of analyzing multi-scale data, we provided mrCOSTS, a robust, unsupervised method capable of diagnosing these data with minimal hyperparameter tuning. It adapts the mrDMD framework for the complexities of real multi-scale data.  In exceptionally challenging real-world data considered here, it already shows its value as an unsupervised algorithm by extracting hitherto unknown spatio-temporal patterns in well studied systems.

Specifically, we demonstrated mrCOSTS on data from a range of fields ranging from biology (neurology) to climate/weather, which were explicitly chosen for their well-known multi-scale dynamics. It was used to successfully and trivially diagnose these complex multi-scale data. For each system we recovered known dynamics and accurately reconstructed the data, while also revealing features which could not be previously characterized due to the complexities of the underlying dynamics.  

In the SST example, we diagnosed a high-dimensional signal with processes spanning orders of magnitude with non-stationarity properties. The known time dynamics were recovered along with their coherent spatial patterns, a useful result on its own. We additionally discovered three time bands outside the normally defined ENSO temporal scales, which were critical for explaining the time series of the ENSO anomaly as well as the extreme behavior of the 2015-2016 ENSO event (Fig. \ref{fig:ENSO-reconstruction}). In the neurology example we deconstructed the LFP signal into coherent spatiotemporal patterns of well-known frequency bands. We revealed for the first time the translating properties (Fig. \ref{fig:monkey}; Movie S1). In the MoBL example, we uncovered dynamics that were previously unknown, revealing a complex interplay of processes that challenged our intuitive expectations (Fig. \ref{fig:MoBL-reconstruction}; Movie S2). In each of these cases, either no or exceptionally limited hyperparameter tuning was performed to highlight the unsupervised diagnostic power provided by mrCOSTS. Of particular note, noise and instrument errors did not hinder analysis for any of the systems.

Diagnosing coherent spatiotemporal patterns with mrCOSTS has untapped benefits, as we only demonstrate the general application of the method here. For example, in the SST example, the 11.1 year band covers a wide range of latitudes around the equator and extends well beyond the Ni\~no 3.4 box (Fig. \ref{fig:ENSO-reconstruction}d), while the 1.4 year band is tightly bound to the equator and extends from the western edge of Ni\~no 3.4 box to the eastern Pacific (Fig. \ref{fig:ENSO-reconstruction}i). The decomposed spatiotemporal patterns should enable a more robust description of ENSO events compared to the wide range of metrics currently needed to describe the complex spatiotemporal patterns of this system \cite{trenberthIndicesNinoEvolution2001, santosoDefiningCharacteristicsENSO2017, lheureuxObservingPredicting20152017, adamsonSituatingNinoCritical2022}.

A robust scale separation for multi-scale data has been a long sought after goal in a wide range of fields from atmospheric turbulence \citep{Sun2015}, weather \citep{honnertAtmosphericBoundaryLayer2020}, ocean dynamics \citep{xiaoReconstructionSurfaceKinematics2023, wangDynamicalDecompositionMultiscale2023}, climate \citep{schurerSeparatingForcedChaotic2013, frankcombeSeparatingInternalVariability2015, Deser2016a}, neurology \citep{zaleskyTimeresolvedRestingstateBrain2014, pedersenSpontaneousBrainNetwork2017, sorrentinoDynamicalInteractionsReconfigure2023}, and others not explored here. Consequently, mrCOSTS has substantial value over other methods for its scale separation properties alone. The bands, $\mathcal{G}_p$, enable analysis requiring isolating processes describable by specific dynamics \citep[similar to][]{wangDynamicalDecompositionMultiscale2023}. The method has additional potential for enabling the retrieval of physics that emerge from complex interactions across scales. For instance, the interpretable DMD model mrCOSTS is based on can facilitate finding separate models capable of describing the different scales and the coupling between them. 

Multi-scale data are a major challenge in many disciplines, especially those defining the Grand Challenges in Science and Engineering \citep{omennGrandChallengesGreat2006}. In order to drive hypothesis testing and process understanding forward we must be able to discover the patterns and processes within these complex data. Recovering previously unknown patterns and complex scale interactions is a task well-suited to analysis by data-driven modeling. But, the complexities of these data have somewhat blunted the potential of such approaches. Providing a robust, unsupervised method of scale separation and diagnoses moves us substantially closer to the goal of using data to derive the complex, emergent laws of multi-scale systems, but does not solve the problem outright. Substantially more work must be conducted to determine how to best use and improve these diagnoses.

\section{Methods}\label{sect:methods}
The python implementation of mrCOSTS is available as part of the PyDMD package \citep{ichinagaPyDMDPythonPackage2024} and is one of the major updates comprising version 1.0 of the package. The implementation includes plotting routines, parameter sweeps, in/out operations for off-line analysis of the fitted model, and tutorials on real and toy data. PyDMD includes automated testing, an active user community, and support for a wide range of other DMD models as well, providing a rich framework for analysis.

\subsubsection*{Practical considerations}
While mrCOSTS's ability to diagnose multi-scale physics has the demonstrated potential to unlock many difficult to analyze types of data, several notes of caution should be observed. All other DMD variants naturally provide an equation-free method for prediction. mrCOSTS does not provide this ability and further the method loses reliability near the the edges of the time domain \cite{Dylewsky2019_dmd}.

Unlike other decomposition like MRA, mrCOSTS does not perfectly decompose the input signal. This has advantages, for instance the exclusion of white noise from fitting. However, some dynamics may not be fit at all so the fit must be carefully evaluated (as in Fig. S1 and S2). Similarly, mrCOSTS does not yet have a method for significance testing. Instead, one must assess which bands are relevant for particular subregions or processes as was done via visual inspection for diagnosing ENSO.

\subsubsection*{Distinction from other methods}

Previous attempts at diagnosing multi-scale data using data-driven modeling required substantial hand tuning and expert intervention \citep[e.g.][]{Kutz_DMD-textoobk_2016, Dylewsky2019_dmd}. Initial developments were done on synthetic data on which we knew the underlying dynamics. Applying them to new and complex data was found to be insufficient, although efforts have been made to account for intermittent or transient phenomenon as well as non-stationary data~\cite{ferreNonStationaryDynamicMode2023}.  Despite such efforts, the need for a unified framework to handle the multitude of difficulties in multi-scale data remained. Indeed, the original sliding mrDMD algorithm performed well for toy model data with regular oscillators and well-separated scales. However, the algorithm struggled when attempting to decompose the types of multi-scale data presented in this manuscript. Many of the overarching concepts from \cite{Kutz2016_dmd} and \cite{Dylewsky2019_dmd} were critical starting points and the principles introduced in these works still remain, on the whole, valid and relevant to mrCOSTS. The algorithmic changes necessary for mrCOSTS to provide the diagnosis of real multi-scale data highlight the importance of testing algorithms on real data and not just synthetic or simple examples.

Another innovation for multi-scale data, non-stationary DMD, was not constructed to separate the various time scales \cite{ferreNonStationaryDynamicMode2023}. The improvements through mrCOSTS created an autonomous and objective retrieval requiring little hyperparameter tuning and expert intervention. As such, this type of unsupervised diagnostic is unprecedented.  

\subsubsection*{Algorithm details}

Separating frequency bands requires transforming $\omega$. Originally in \citep{Dylewsky2019_dmd} clustering was performed only on $\omega\omega^*$ (where we use $\omega^*$ to denote the complex conjugate), which includes the contribution of both $Re(\omega)$ and $Im(\omega)$. For small amplitude modes, the real and imaginary components can alternate between adjacent windows with similar values (e.g., $10 + 0i$ vs $0 + 10i$). The transformation $\omega\omega^*$ treats these values as identical, contaminating the frequency band identification. mrCOSTS supports multiple methods of transforming $\omega$ for clustering frequency bands: $|Im(\omega)|$, $|Im(\omega)|^2$, and $log_{10}(|Im(\omega)|)$. We find $|Im(\omega)|$ is often the appropriate choice of transformation for the local scale separation. For the global frequency band identification, clustering must be performed using the $log_{10}(|Im(\omega)|)$ transformation. Otherwise, the frequency bands, which can span order of magnitudes of scales, cannot be reliably separated.

The DMD model fit to each window is the variable projection optimized DMD \cite{Askham2018}. We introduce a powerful innovation which enables placing arbitrary constraints on the eigenvalue solutions and is now implemented in PyDMD \cite{ichinagaPyDMDPythonPackage2024}. We found in particular that requiring the real component of the eigenvalues to be small, but not exactly zero, enables robust fitting of data. Most physical systems, while the data may be non-stationary or non-linear, do not have unconstrained growth or decay and by imposing this eigenvalue constraint we recognize this feature. Without this constraint, or alternatively by imposing $Re(\omega)=0$, individual windowed DMD fits can be poor, especially when fitting the smaller decomposition levels.

We also make note of the strategy for finding initial values for the eigenvalues. Initial values for eigenvalues are not shared between windows. Instead, each window is independently initialized. This strategy leads to better overall fits across all windows at the cost of computation time. Other initialization methods regularly failed for data without well-separated frequency bands, as is commonly the case for real multi-scale data. Using a fixed set of initial eigenvalues for all windows functions for systems with regular oscillators common in toy model data but is not generally recommended.

\subsubsection*{Hyperparameters}
 A wide variety of data can be fit with little to no hyperparamter tuning. We explicitly refrained from hyperparameter tuning, except for adding a decomposition level in the neurology and MoBL examples, to demonstrate this point. The most sensitive parameters are the svd rank of the DMD fit ($r$), the window size for each decomposition level ($w$) as well as the number of decomposition levels ($N$), and the constraints on the eigenvalues. When encountering a bad fit, often manifested as a small number of individual windows with poor fits, these are the hyperparameters to adjust.

 Objective measures for deriving an optimal $r$ based on energy arguments dramatically lower the performance of mrCOSTS relative to a fixed $r$ in a given decomposition level. Instead the rank should be fixed as an even number. The number of resolved bands, $p$, generally scales with $\frac{r}{2}$.

The window size for each decomposition level determines the largest resolvable high frequency component. e.g., a small window size can fit higher frequency components while a larger window size can fit lower frequency components \cite[similar to][]{rahamanSpectralBiasNeural2019}. Dyadic scaling of the window size is not necessary (see the MoBL for example) and qualitatively similar results can be recovered for varying window sizes. Since strict scaling is not necessary, window sizes can be chosen to target scales of interest. Typically, the smallest window size is the most finicky since the window length needs to be slightly larger than $r$ but setting it to be too large means omitting the finest scale features. Generally, the windowed DMD fit scales closer to the window size and omits the highest frequency components creating a trade-off between a large enough $r$ and a small enough $w$ to capture the highest frequency components. Missing frequencies in $\mathcal{G}_p$ can usually be solved by increasing $r$ for the relevant decomposition level(s). 

The number of frequency bands, both for the global and local separations, can be provided as an \textit{a priori} expectation or found objectively using k-means clustering from scikit-learn \citep{scikit-learn2011} with a hyperparameter sweep. The optimal number of frequency bands is found using the Silhouette score \cite{rousseeuwSilhouettesGraphicalAid1987}, which is well-suited for clustering of the flattened $\omega$ arrays. Other metrics and clustering algorithms were tested and found to perform substantially worse.

\subsection*{Data}
All data are publicly available.

\subsubsection*{Sea Surface Temperatures}

We use the Hadley Centre Sea Ice and Sea Surface Temperature data set \cite{raynerGlobalAnalysesSea2003} using the SST data between 30$^{\circ}$ S and 30$^{\circ}$ N and 100$^{\circ}$ W and 120$^{\circ}$ E. HadISST data \cite{raynerGlobalAnalysesSea2003} are publicly available at https://www.metoffice.gov.uk/hadobs/hadisst/. We calculate SST$^*$ relative to a moving 30 year historical period and smoothed by a rolling 3 month average \citep[as in][]{noaacliamtepredictioncenterNOAAClimatePrediction}. When evaluating within the Ni\~no 3.4 box, SST$^*$ is relative to the spatial mean within the box \citep[consistent with][]{noaacliamtepredictioncenterNOAAClimatePrediction, trenberthkevinNinoSSTIndices}, while for the entire domain the anomaly is calculated relative to each grid point individually.

For the decomposition, no hyperparameter tuning was performed. We specified the rank of the decomposition for each level ($r$=8), decomposition windows of [16, 32, 64, 128, 256, 512, 1024] (dyadic scaling), and a slide of approximately 10\% between windows. The initial window length was selected to be larger than the decomposition rank but short enough that seasonal and subseasonal modes could be identified.

The total error of the reconstruction over the entire period was 11\%, with most of the error located at the sub-seasonal scales and edge effects (Fig. S1). The fit could be improved by increasing $r$, especially for the smaller windows, and constraining the eigenvalues to have a small $Re(\omega)$. This was not performed as it was unclear what would constitute a small $Re(\omega)$ until after the first decomposition was performed and due to the unknown trade-off between $r$ and the window length of the first decomposition level.

\subsubsection*{Neurology}

The neurology data can be accessed through the DOI 10.12751/g-node.f83565 \cite{riehleMappingSpatiotemporalStructure2013, brochierMassivelyParallelRecordings2018}. The data for trail 2 of monkey L were quality-controlled, removing points with exceptionally large or small temporal deviations or uncharacteristic jumps, resulting in 7 of the 96 electrodes being discarded and a total of 11 missing points in the spatial map of the electrodes (Fig. \ref{fig:monkey}d-g).

The decomposition was performed from the start of the trial to the end of the reward being dispensed (duration about 3.6 s). The decomposition window sizes were [50, 100, 250, 500, 1000, 2000] microseconds with a slide of 10\% of the window length. Each window length was fit using $r$ = 8 and we forced each level to find 4 frequency bands for the local scale separation and objectively found the number of global bands using the built-in hyperparameter sweep from PyDMD (i.e., $p$ was not set \textit{a priori}).

The mrCOSTS reconstruction had an error of 21\% for the full period and 14\% when excluding the edge effects (approximately 0.15 s at the start and end). Most of the error was a consequence of dropping white noise, poorly resolved high frequencies (Fig. \ref{fig:monkey}d), and edge effects.

\subsubsection*{Mountain Boundary Layer}

The LIDAR co-planar data can be accessed through the Zenodo repository at the DOI 10.5281/zenodo.7212801 \cite{pfisterTEAMxPC22Alpine2024}. The co-planar retrieval \cite{Kalthoff2014} was conducted using two LIDARS placed on the north valley slope, opposite of the tributary valley in-flow. The LIDAR retrieval yielded only radial velocity relative to the instrument, but using co-located observations, the $u$ and $v$ velocity components were retrieved. A complete observation of the co-planar area, i.e. the data time step, was made every 2 min. Gaps for individual grid points less than 4 observations long (maximum duration of 6 minutes) were filled using a linear interpolation in time.

We decomposed a single night from 2000 the night before to 0715 the next morning. The window lengths were [24, 48, 96, 192, 360, 480] minutes (the longest window was 8 hours in length) with the slide between windows being 10\% of the window length. Each window length was fit using $r$ = 8 and we forced each level to find 4 frequency bands for the local scale separation and objectively found the number of global bands using the built-in hyperparameter sweep from PyDMD (i.e., $p$ was not set \textit{a priori}).

The mrCOSTS fit had an error of 23\% and 26\% for the $u$ and $v$ components, respectively. The majority of the error was located at the beginning and end of the decomposition due to the edge effects. The error decreases to 13\% and 15\% when excluding the first and last twenty minutes (Fig S2). The mrCOSTS decomposition recovered the complex temporal dynamics except at time scales smaller than 8 min (Fig. S2a-b) as a result of the relatively coarse observation time step, instrument noise, observation gaps, and prioritizing larger scales with the hyperparameters.

The decomposition accurately recovered $u$ and $v$ low-order statistical moments (Fig. S2c-f), demonstrating how mrCOSTS can robustly decompose multiple variables through their shared covariance. As a result, $<U>$ was also accurately reconstructed from the mrCOSTS fit (Fig. \ref{fig:MoBL-reconstruction}a). The error in the mrCOSTS fit was mostly caused by missing the small time scales, removing white noise, and the edge effects.

\section*{Acknowledgements}
KL was funded by the Austrian Science Fund (FWF) [10.55776/ESP214]. For open access purposes, the author has applied a CC BY public copyright license to any author-accepted manuscript version arising from this submission.  SMI and JNK acknowledges support from the National Science Foundation AI Institute in Dynamic Systems (grant number 2112085).

\bibliographystyle{naturemag}
\bibliography{arxiv-library}

\end{document}